\documentclass[12pt, dvips, twoside]{amsart}

\swapnumbers\newtheorem{theorem}{Theorem}[section]

\newtheorem{corollary}[theorem]{Corollary}

\newtheorem{proposition}[theorem]{Proposition}

\newtheorem{lemma}[theorem]{Lemma}

\newtheorem{definition}[theorem]{Definition}

\newtheorem{remark}[theorem]{Remark}
\newtheorem{remarks}[theorem]{Remarks}

\newtheorem{claim}[theorem]{Claim}

\numberwithin{equation}{section}

\usepackage{amssymb,latexsym, amsmath, amscd, graphicx, epsfig, psfrag}

\def\po{\parindent 0pt}
\def\p{{\po \it \s Proof. }}
\def\pp{{\po \it \m Proof}}
\def\s{\smallskip}
\def\m{\medskip}
\def\ni{\noindent}
\def\ra{\rightarrow}

\def\cat{\operatorname{cat}}
\def\ccat{\overline{\operatorname{cat}}\,\!}
\def\gcat{\cg \text{-} \!\operatorname{cat}} 
\def\1cat{\cg_1 \text{-} \operatorname{cat}} 
\def\2cat{\cg_2 \text{-} \operatorname{cat}} 
\def\gccat{\cg \text{-}\,\! \overline{\operatorname{cat}}\,\!} 

\def\Crit{\operatorname{Crit}}
\def\mmod{\operatorname{mod}}

\def\RR{\mathbb {R}}
\def\NN{\mathbb {N}}
\def\CC{\mathbb {C}}
\def\ZZ{\mathbb{Z}}

\def\ga{\alpha}
\def\gb{\beta}
\def\gg{\gamma}
\def\gd{\delta}
\def\eps{\varepsilon}

\def\gf{\varphi}

\def\gs{\sigma}

\def\BB{{\mathcal B}}
\def\cg{{\mathcal G}}

\newcommand{\proofend}{\hspace*{\fill} $\Box$\\}
\newcommand{\diam}{\hspace*{\fill} $\Diamond$}

\begin{document}

\title[Lusternik--Schnirelmann theory on general spaces.]{
Lusternik--Schnirelmann theory on general spaces}

\date{\today}

\author{Yu.\ B.\ Rudyak}
\address{(Yu.\ B.\ Rudyak)  Max--Planck Institut f\"ur Mathematik, Vivatsgasse 7, 53111 Bonn, Germany}
\email{rudyak@mpim-bonn.mpg.de}
\author{F.\ Schlenk}
\address{(F.\ Schlenk) ETH Z\"urich, CH-8092 Z\"urich, Switzerland}
\email{felix@math.ethz.ch}

\begin{abstract}
We extend Lusternik--Schnirelmann theory to pairs $(f, \gf)$, where
$\gf$ is a homotopy equivalence of a space $X$, 
$f$ is a function on $X$ which decreases along $\gf$ and
$(f, \gf)$ satisfies a discrete analog of the Palais--Smale condition.
The theory is carried out in an equivariant setting.
\end{abstract}

\maketitle

\section{Introduction}

\ni
The Lusternik--Schnirelmann category $\cat_X \!A$ of a subset $A$ of a
topological space $X$ is the minimal number of open and in $X$
contractible sets which cover $A$ \cite{LS, F}. 
We set $\cat X = \cat_X \! X$.
Lusternik and Schnirelmann \cite{LS} 
proved that 
for any smooth function $f$ on a closed manifold $M$
\begin{equation*}
 \sum_{d \,\in\, \RR} \cat_M  K \cap f^{-1}(d) \,\ge\, \cat M 
\tag{a}
\end{equation*}
where $K$ denotes the set of critical points of $f$.
In particular, $f$ has at least $\cat M$ critical points, and if $M$ is
connected, then                
\begin{equation*}
\text{$K$ is infinite or $f$ has at least $\cat M$ critical levels.}
\tag{b}
\end{equation*}
For generic $f$, (a) can be written as
\begin{equation*}
 \cat K \,\geq\, \cat M .
\tag{c}
\end{equation*}
In \cite{P4}, the estimates (a) and (b) have been generalized to a class
of Banach manifolds $M$ and functions $f \colon M \ra \RR$ which satisfy
the Palais--Smale condition,
and more recently, versions for continuous functions on certain metric
spaces have been found (see \cite{CD}).
In this paper, we extend (a), (b) and (c) to more general
situations.

\s
Consider a topological space $X$, a continuous map
$\gf \colon X \ra X$ and a continuous function $f \colon X \ra \RR$.
We say that $f$ is a {\it Lyapunov function}\, for
$\gf$ if $f(\gf (x)) < f(x)$ whenever $x$ is not a fixed point of $\gf$.

\begin{definition}
{\rm Let $f$ and $\gf$ be as above and let $Y$ be a subspace of $X$.
We say that the pair $(f, \gf)$ satisfies condition (D) on $Y$ if $f$
is a Lyapunov function for $\gf$ and the following holds.
\begin{quote}
If $A$ is a subset of $Y$ on which $|f|$ is bounded but on which
$f(y) - f(\gf(y))$ is not bounded away from zero, then there is a fixed
point of $\gf$ in the closure of $A$.
\end{quote}
}
\end{definition}
\s
\ni
Observe that if $Y$ is compact or if $f$ is proper then the pair
$(f,\gf)$ satisfies condition (D) whenever $f$ is a Lyapunov function for
$\gf$.
Condition (D) is a discrete analog of the Palais--Smale condition for
$C^1$ functions on Banach manifolds. In fact, we show in Proposition
\ref{p:app} that in the case of a Hilbert manifold, condition (D)
generalizes the Palais--Smale condition.

A topological space $X$ is called {\it weakly locally contractible}\, 
if each point of $X$ has a neighborhood which is contractible in $X$.
Observe that this property is equivalent to $\cat X$ being definite.
A space $X$ is called {\it binormal}\, if $X \times [0,1]$
(and hence $X$) is normal. In particular, every metric space and every
compact space is binormal.
A space $X$ is called an {\it absolute neighborhood retract}\, (ANR) if
for each closed subset $Z$ of a normal space $Y$ every map $\rho \colon Z
\to X$ admits an extension to a neighborhood of $Z$.
CW-complexes and topological manifolds are ANR's.
Every binormal ANR is weakly locally contractible 
(see \text{Lemma \ref{l4:2}}).

\m
If some category $\cat_X A$ is not finite, i.e., infinite or indefinite, 
we write $\cat_X A = \infty$. 
We agree that $\infty \ge \infty$ and that $\infty \ge n$ for every $n
\in \ZZ$.
Given a real valued function $f$ on a space $X$, we set
\[
f^a = \{ x \in X \,|\, f(x) \le a \}.
\]

\begin{theorem}  \label{t0}
Consider a space $X$ and a homotopy equivalence $\gf$ of
$X$. Let $F$ be the set of fixed points of $\gf$. Assume
that the function $f \colon X \ra \RR$ is bounded below and that $(f,
\gf)$ satisfies condition {\rm (D)} on $f^b$ for some $b \in \RR$.

\s
\begin{itemize}
\item[(a)]
We always have
$$
\sum_{d \,\le\, b} \cat_X  F \cap f^{-1}(d)
\, \ge \, \cat_X f^b.
$$
\item[(b)]
If $X$ is connected and weakly locally contractible, then $F \cap f^b$
is infinite or $f(F) \, \cap \, ]- \infty, b]$ contains at least $\cat_X
f^b$ elements. 

\s
\item[(c)]
If $X$ is a binormal ANR and $f(F)$ is discrete, then
$$
\cat F \cap f^b \,\ge\, \cat_X f^b.
$$
\end{itemize}
If $\gf$ is homotopic to the identity map, then these
estimates also hold for $b = \infty$ (in which case $f^b = X$).
\end{theorem}

Our main results generalize Theorem \ref{t0} in two directions.
Firstly, we overcome the assumption of $f$ being bounded below: 
Given $- \infty < a < b < \infty$, Theorem \ref{t5:1} provides 
a lower bound for the number and the category of the rest
points of a homotopy equivalence $\gf$ in the slice 
$\{ x \in X \mid a < f(x) \le b \}$,
and in Theorem \ref{t6:1} we show that if $\gf$ is homotopic to the
identity, then the categorical estimates in Theorem \ref{t5:1} can be
replaced by various finer ones, which also hold for $b = \infty$.
Secondly, all results readily extend to equivariant situations where
a compact Lie group acts on $X$.

Finally, we observe in Section 8 that if $\gf$ is a homeomorphism or the
time-1-map of a flow, then some of the previously established estimates can be
improved by using different invariants of Lusternik--Schnirelmann type.

\m
Notice that even if $\gf$ is homotopic to the identity, 
we in general do not require that there exists a homotopy between the
identity and $\gf$ which decreases along $f$. 
Therefore, our theory should also prove useful in understanding the
fixed point set of other than gradient-like systems such as those
arising in simulated annealing. 

The need of extending Lusternik--Schnirelmann theory to spaces more
general than manifolds was demonstrated by various proofs of the
Arnold conjecture about the number of fixed points of Hamiltonian
symplectomorphisms \cite{H, HZ, R}. In these proofs different variants
of the Lusternik--Schnirelmann category were considered. We thus
develop
an axiomatically defined version of Lusternik--Schnirelmann theory
(Theorem \ref{t2}) and derive Theorem \ref{t5:1} and Theorem \ref{t6:1} 
from Theorem \ref{t2}.
In \cite{RS} this axiomatic approach will be used to develop a
Lusternik--Schnirelmann theory for certain symplectic and contact
manifolds.

\m All spaces are assumed to be Hausdorff, and all maps and functions are
assumed to be continuous.

\section{An axiomatic version of LS-theory}

\begin{definition}
[\rm cf.\ \cite{LS, CP1, CP2}]
\rm
An {\it index function}\, on a topological space $X$ is a function
$ \nu \colon 2^X \times 2^X\to \NN \cup \{0\}
$
satisfying the following axioms:

\s
{\bf (monotonicity)}\, If $A \subset B$, then $\nu (A,Y) \leq \nu
(B,Y)$ for all $Y\subset X$;

\s
{\bf (continuity)}\, For every closed subset $A \subset X$ there
exists a neighborhood $U$ of $A$ such that $\nu (A,Y) = \nu (U,Y)$ for
all $Y \subset X$;

\s
{\bf (mixed subadditivity)}\, For all $A, B, Y \subset X$,
\[
\nu(A\cup B,Y) \,\leq\, \nu (A,Y) + \nu (B,\emptyset).
\]

\m
\ni
Given a map $\gf: X \to X$ and a subset $Z$ of $X$, an index function $\nu$ is {\it
$(\gf, Z)$-supervariant}\, if

\s
{\bf ($(\gf, Z)$-supervariance)}\, $\nu (\gf (A),Z) \ge \nu (A,Z)$ for all $A
\subset X$.

\end{definition}

Below we write $\nu(A)$ instead of $\nu(A, \emptyset)$.

\begin{lemma}
[\rm cf.\ \cite{LS}]
\label{l2:1}
Consider a self-map $\gf \colon X \to X$ on the space $X$. Denote the
set of fixed points of $\gf$ by $F$. Let $f \colon X \ra \RR$ be a
function and let $a < b$ be real numbers such that $(f,\gf)$ satisfies
condition {\rm (D)} on $f^{-1}[a,b]$. Then the following assertions
hold.

\s
{\rm (i)} Suppose that $f^{-1} [a,b[ \, \cap \, F$ is empty.
Then, given any neighborhood $U$ of $f^{-1} (b) \cap F$, there exists $n \in
\NN$ such that $\gf^n (f^b \setminus U) \subset f^a$. In particular,
$\nu(f^b \setminus U, Z) \leq \nu (f^a, Z)$ for every $(\gf, Z)$-supervariant index function $\nu$, and hence $\nu(f^b \setminus U, Z) = \nu (f^a, Z)$
whenever $f^a \cap U=\emptyset$.

\s
{\rm (ii)} Suppose that there is $\eps>0$ such that
$f^{-1} ]a,a+ \eps[ \, \cap \, F$ is empty.
Then, given any neighborhood $U$ \!of $f^a$, there exists $\gd>0$ such
that $\gf (f^{a+\gd}) \subset U$.
In particular, for every $(\gf,Z)$-supervariant index function $\nu$,
there exists $\gd>0$ such that $\nu(f^{a+\gd}, Z) = \nu(f^a, Z)$.
\end{lemma}

\p (i)
Since $(f, \gf)$ satisfies condition (D) on $f^{-1}[a,b]$ and since
$\gf$ has no fixed point on the closed set $f^{-1}[a,b]
\setminus U$, there exists $\delta_1 >0$ such that
\begin{align*}
f(x) - f(\gf (x)) & \geq \delta_1
        \quad \text{for every } x \in f^{-1}[a, b] \setminus U. \\
\intertext{If $b - \gd_1 \le a$, we are done. Otherwise, let $\gd_2 >0$
be
such that}
f(x) - f(\gf (x)) & \geq \delta_2  \quad \text{for every } x \in
f^{-1}[a, b-\gd_1] .
\end{align*}
Choose $n \in \NN$ so large that $(n-1) \gd_2 \ge b-\gd_1-a$.
Then $\gf^{n} (f^b \setminus U) \subset f^a$.
Indeed, assume that $f(\gf^{n}(x)) >a$ for some $x \in f^b \setminus
U$. Then $\gf^k (x)
\in f^{-1} [a, b-\gd_1]$ for $k = 1, \dots, n$. Hence,
\begin{eqnarray*}
f(\gf^{n} (x)) &=& f(\gf^1(x)) + \sum_{k=1}^n f(\gf^{k+1} (x)) -
                        f(\gf^k (x))  \\
                 &\le & b-\gd_1 - (n-1) \gd_2 \, \le \, a,
\end{eqnarray*}
a contradiction.

Assume that $\nu$ is a $(\gf,Z)$-supervariant index function. Choose $n$ such
that $\gf^{n} (f^b \setminus U) \subset f^a$. Then, by
$(\gf,Z)$-supervariance and monotonicity,
\[
\nu (f^b \setminus U,Z) \le \nu (\gf^n (f^b \setminus U),Z) \le \nu (f^a,Z).
\]

\s
(ii) Fix a neighborhood $U$ of $f^a$ and
choose a sequence of real numbers
\[
a + \eps > a_1 > a_2 > \cdots > a_n > \cdots
\]
which converges to $a$.
We claim that $\gf (f^{a_n}) \subset U$ for some $n$. Arguing by
contradiction, suppose that
$
A_n = f^{a_n} \setminus \gf^{-1} (U) \neq \emptyset$\, for all $n $.
Then $f(x) - f(\gf (x)) \le a_n - a$\, for all $x \in A_n$, and so
$\inf
\{f(x) - f(\gf (x)) \mid x \in A_1 \} =0$. Since $A_1$ is closed,
condition (D) implies that $A_1 \cap F \neq \emptyset$. This
contradicts $f^{-1} ]a,a+\eps[ \, \cap \, F = \emptyset$. So, $\gf
(f^{a_n}) \subset U$ for some $n$. Set $\gd = a_n - a$.

Assume that $\nu$ is a $(\gf,Z)$-supervariant index function.
By continuity, there exists a neighborhood $U$ of $f^a$
such that $\nu (f^a,Z) = \nu (U,Z)$. According to what we said above, there
exists $\gd>0$ such that $\gf (f^{a+\gd}) \subset U$. Then
\[
\nu (f^{a}, Z) = \nu (U, Z) \ge \nu (\gf (f^{a+\gd}), Z) \ge \nu (f^{a+\gd}, Z)
\ge \nu (f^{a}, Z)
\]
and so $\nu (f^{a+\gd}, Z) = \nu (f^a, Z)$.
\proofend

\ni
The Lusternik--Schnirelmann Theorem, from the modern standpoint, is

\begin{theorem}  \label{t2}
Consider a self-map $\gf \colon X \to X$on the topological space $X$.
Let $F$ be the set of fixed
points of $\gf$ and let $f$ be a Lyapunov function for $\gf$. Let $a <
b$ be real numbers such that $(f,
\gf)$ satisfies condition {\rm (D)} on $f^{-1}[a, b]$.
Suppose that
$f(F) \, \cap \, ]a,b] = \{ d_1, \dots, d_m \}$ is a finite set.
Then
$$
\nu (f^a, f^a) \,+\,
\sum_{i=1}^m \nu \left( F \cap f^{-1}(d_i) \right) \,\geq\, \nu(f^b, f^a).
$$
for any $(\gf, f^a)$-supervariant index function $\nu$.
\end{theorem}

\p
Denote $\nu(A,f^a)$ by $\mu(A)$.
If $\mu (f^b) = \mu (f^a)$, there is nothing to
prove. So assume that $\mu (f^a) < \mu (f^b)$. Set
\[
c_k = \inf \{ c \in \RR \bigm| \mu(f^c) \ge k \}, \qquad
k = \mu(f^a)+1, \ldots, \mu(f^b).
\]
By monotonicity of $\mu$,
\begin{equation}  \label{e1}
a \le c_{\mu(f^a)+1} \le \dots \le c_{\mu(f^b)} \le b .
\end{equation}
\begin{claim}  \label{c:k}
For all $k = \mu(f^a)+1, \ldots, \mu(f^b)$,
 \begin{eqnarray}
   \mu(f^{c_k-\eps})  &<&       k   \quad \text { for all \,$\eps > 0$},
                                         \label{e:nu1}  \\
   \mu(f^{c_k})       &\ge&     k .      \label{e:nu2}
 \end{eqnarray}
\end{claim}

\p
Inequality \eqref{e:nu1} follows from the definition of $c_k$.

We prove \eqref{e:nu2}. By \eqref{e1}, $c_k\leq b$. Suppose first that
$c_k <b$.
Since $f(F) \, \cap \, ]a,b]$ is finite, we find $\eps >0$ such that
$f(F) \, \cap \, ]c_k, c_k+\eps[ \, = \emptyset$, i.e., $f^{-1}
]c_k,c_k+\eps[ \, \cap \, F = \emptyset$. By Lemma \ref{l2:1} (ii) there
exists $\gd >0$ such that $\mu (f^{c_k+\gd}) = \mu (f^{c_k})$. But, by
monotonicity, $\mu (f^{c_k+\gd}) \ge k$, and so $\mu (f^{c_k}) \ge k$.
Suppose now $c_k =b$. Then $\mu (f^{c_k})=\mu (f^b) \geq k$ as well.
\diam

\m
Plugging $k = \mu (f^a)+1$ into \eqref{e:nu2}, we find $\mu ( f^{c_{\mu
(f^a) +1}}) \ge \mu (f^a)+1$. Hence,
\begin{eqnarray}  \label{e:a}
a<c_{\mu (f^a)+1} .
\end{eqnarray}

\begin{claim}  \label{c:point}
$f^{-1}(c_k)$ contains at least one point of $F$.
\end{claim}

\proof
Suppose the contrary. In view of \eqref{e1}, \eqref{e:a} and
the finiteness of $f(F) \, \cap \, ]a,b]$, we
then find $\eps >0$ such that $f^{-1}[c_k
-\eps, c_k]$ does not contain points of $F$. By
Lemma \ref{l2:1}\,(i), there exists $n$ such that $\gf^n (f^{c_k})
\subset f^{c_k -\eps}$. Hence, by $(\gf,f^a)$-supervariance, monotonicity and
\eqref{e:nu1},
\[
\mu(f^{c_k}) \le \mu (\gf^n (f^{c_k})) \le \mu(f^{c_k-\eps}) < k.
\]
This contradicts \eqref{e:nu2}.
\diam

\begin{lemma}
[{\rm cf.\ \cite[II, \S \!\! 4]{LS}, \cite[Lemma 19.12]{DFN}}]
\label{l2:2}
If $c_k = c_{k+1} = \cdots = c_{k+r}$, then
\[
\nu (F \cap f^{-1}(c_k)) \geq r+1.
\]
\end{lemma}

\p
Set $A = F \cap f^{-1}(c_k)$. By continuity, there is a neighborhood
$U$ of $A$ with $\nu(U)=\nu(A)$. Arguing as above we find $\eps >0$
such that $f^{-1} [c_k -\eps, c_k[ \; \cap \, F$ is empty. By
Lemma \ref{l2:1}\,(i) there exists $n$ such that $\gf^n (f^{c_k}
\setminus U) \subset f^{c_k -\eps}$. So, by $(\gf,f^a)$-supervariance, monotonicity
and (\ref{e:nu1}), $\mu(f^{c_k} \setminus U) < k$. Now,
by monotonicity, subadditivity and \eqref{e:nu2},
\begin{eqnarray*}
\mu(f^{c_k} \setminus U) + \nu (U) &=&
\mu (f^{c_{k+r}} \setminus U) + \nu(U) \\
 &\ge& \mu (f^{c_{k+r}} \setminus U) +
                 \nu(f^{c_{k+r}} \cap U) \\
 &\ge& \mu (f^{c_{k+r}}) \\
 &\ge& k+r,
\end{eqnarray*}
and so $\nu(U)>r$. Thus, $\nu(A) \geq r+1$.
\diam

\m We continue the  proof of the theorem. We have
\begin{eqnarray*}
c_{\mu(f^a)+1} \,= \cdots =\, c_{i_1} &<&  c_{i_1+1} \,= \cdots =\,
c_{i_2} \\
       &<& \cdots <\, c_{i_n+1} \,= \cdots =\, c_{i_{n+1}} \,=\,
c_{\mu(f^b)}.
\end{eqnarray*}
Set $F_k = F \cap f^{-1} (c_{i_k})$ and $i_0 = \mu(f^a)$.
Then, by Lemma \ref{l2:2}, $\nu(F_k) \ge i_k-i_{k-1}$. Thus
$$
\sum_{k=1}^{n+1}\nu(F_k) \ge i_{n+1} - i_0 = \mu (f^b) -\mu (f^a).
$$
In view of \eqref{e:a} and Claim \ref{c:point},  $\{
c_{\mu(f^a)+1},\dots, c_{\mu (f^b)}\} \subset \{ d_1,\ldots,
d_m \}$, and so Theorem \ref{t2} follows.
\proofend

\section{Relative equivariant category}

\ni
Let $G$ be a compact Lie group. A {\it $G$-space}\, is a space $X$
together with a continuous action $G\times X \ra X$, $(g,x) \mapsto
gx$. 
A {\it subspace}\, of a $G$-space $X$ is a $G$-invariant subspace.
If $X_1$ and $X_2$ are $G$-spaces, then a {\it $G$-map}\,
$\gf \colon X_1 \ra X_2$ is an equivariant map, i.e.,
$\gf (gx) = g \gf (x)$ for all $g \in G$ and $x \in X_1$,
and a {\it $G$-homotopy}\, $H_t \colon X_1 \ra X_2$ is a map
$H \colon X_1 \times I \ra X_2$ such that $H_t \colon X_1 \ra X_2$
is a $G$-map for each $t$.

Let $W$, $Y$ be subspaces of a $G$-space $X$. We say that 
{\it $W$ is $G$-deformable to $Y$}\, 
if there is a $G$-homotopy $H_t \colon W \ra X$ which starts with the 
inclusion and is such that $H_1(W) \subset Y$.
If in addition $H_t (W \cap Y) \subset Y$ for all $t$, we say that
{\it $W$ is $G$-deformable to $Y \mmod Y$}.

\m
Fix a class $\cg$ of homogeneous $G$-spaces,
\[
\cg \subset \{ G/H \mid H \subset G \text{ is a closed subgroup\,\!} \} .
\]
A subspace $A$ of a $G$-space $X$ is called {\it $\cg$-categorical}
if there exist $G$-maps $\ga
\colon A \ra G/H$ and $\gb \colon G/H \ra X$ with $G/H \in \cg$
such that the inclusion $A \hookrightarrow X$ is $G$-homotopic to the
composition $\gb \ga$.

\begin{definition}[\rm cf.\ \cite{Fa, CP2}]
\label{def:gcat}
\rm
Fix a subspace $Y$ of the $G$-space $X$. If $A$ is another subspace of
$X$, we set $\gcat_X (A, Y) =k$ if $A$ can be covered by $k+1$ {\it
open}\, subspaces $A_0, A_1, \dots, A_k$ of $X$ such that

\s
\ni
(i) \quad \;\;\;\; $A_0$ is $G$-deformable to $Y$

\s
\ni
(ii) \quad \;\;\,\,
$A_1, \dots, A_k$ are $\cg$-categorical

\s
\ni
and if $k$ is minimal with this property. If there is no such number
$k$, we set $\gcat_X (A, Y) =\infty$.

The invariant $\gcat_X (A \mmod Y)$ is defined by replacing (i) by

\s
\ni
(i mod) \;\! $A_0$ is $G$-deformable to $Y \mmod Y$, \,and\, $A_0
\supset A \cap Y$.

\m
\ni
We set $\gcat_X A = \gcat_X (A, \emptyset)$ and $\gcat X = \gcat_X X$.
\end{definition}

\begin{remarks}
{\rm
{\bf 1.}
If $G$ acts trivially on $X$ and $\cg$ contains the point $G/G$,
then $A$ is $\cg$-categorical iff it is
contractible in $X$. Therefore, in this case $\gcat_X A$ equals
the classical open Lusternik--Schnirelmann category $\cat_X A$
\cite{LS, B1, F}.

\s
{\bf 2.} If $\cg_1 \subset \cg_2$, then $\1cat_X (A, Y) \ge \2cat_X(A,
Y)$. If $G$ has a fixed point $x$ and $\cg$ does not contain $G/G$,
then $\gcat_X A =\infty$ whenever $x \in A$.

\s
{\bf 3.}
If $G$ acts freely on $X$ and $\cg$ is the full class of homogeneous
$G$-spaces, then $\gcat_X A = \cat_{X/G} A/G$.
In general, however, $\gcat_X A \ge \cat_{X/G} A/G$. E.g., if $\ZZ_2$
acts on $S^1 = \{ z \in \CC \mid |z| =1 \}$ by complex conjugation, then
$\gcat S^1 \ge 2 >1 = \cat S^1 / \ZZ_2$.

\s
{\bf 4.}
Agreeing that $\infty \ge \infty - \infty$ and $0 \ge n - \infty$ for
every $n \in \ZZ$,
we always have
\[
\gcat_X (A \mmod Y) \,\ge\, \gcat_X (A,Y) \,\ge\, \gcat_X A - \gcat_X Y .
\]
If $X = \{ (x,y) \in S^1 \mid y \ge 0 \}$ and $Y = X \cap \{ y =0 \}$,
then 
$$
\cat_X (X \mmod Y) =1>0 = \cat_X (X,Y);
$$ 
and if $X$ is the union of the unit circle $S$ centered at $(0,-1)$ and its translate
centered at $(0,1)$, then $\cat_X (X,S) = 1 > 2-2 = \cat_XX - \cat_X
S$. }
\end{remarks}

\section{Some general topology}

\ni
If $A$ is a subspace of the $G$-space $X$, we set
$\gccat_X A = k$ if $A$ can be covered by $k$ {\it closed}\, categorical
subspaces $A_1, \dots, A_k$ of $X$ and if $k$ is minimal with this
property. If there is no such $k$, we set $\gccat_X A = \infty$. We also
set $\gccat \,\!X = \gccat_X X$, and, if $G$ acts trivially
and $\cg$ contains $G/G$, $\ccat_X A
= \gccat_X A$.

\begin{lemma}  \label{l4:1}
Let $A$ be a closed subspace of a normal $G$-space $X$. Then
\[
\gcat_X \,\!A \, \ge \, \gccat_X A.
\]
\end{lemma}
\p
We may assume that $\gcat_X A$ is finite.
Let $A \subset U_1 \cup \dots \cup U_k$
where each $U_i$ is open and
$\cg$-categorical. Then $\{ X \setminus A, U_1, \dots, U_k \}$ is an
open covering of $X$. According to
\cite[$\S$ \!2, Proposition 20]{AF} there is an open covering $\{
V_0, V_1, \dots, V_k \}$ of $X$ such that $\overline{V_0}
\subset X \setminus A$ and $\overline{V_i} \subset U_i$, $i = 1, \dots,
k$. Hence, $A \subset \overline{V_1} \cup \dots \cup
\overline{V_k}$. By Proposition 1.1.1 of \cite{P1}, every set
$G \overline{V_i} = \{\,\!gv \,|\, g \in G,\, v \in \overline{V_i}\,\}$
is closed. The inclusions $\overline{V_i} \subset  G\overline{V_i}
\subset GU_i \subset U_i$ imply that $A \subset \cup G\overline{V_i}$
and that each $G \overline{V_i}$ is a closed $\cg$-categorical subspace.
Thus, $\gccat_X A \le k$.
\proofend

A $G$-space $X$ is called a {\it $G$-ANR}\, if
for each closed subspace $Z$ of a normal $G$-space $Y$ every $G$-map
$\rho \colon Z \to X$ admits an equivariant extension to a $G$-invariant
neighborhood of $Z$. Notice that $G/H$ is a $G$-ANR for every closed
subgroup $H$ of $G$ \cite[Corollary 1.6.7]{P1}.

\begin{lemma} \label{l4:2}
Let $A$ be a closed $\cg$-categorical subset of a binormal $G$-ANR
$X$. Then $A$ is contained in an open $\cg$-categorical subset of $X$.
\end{lemma}

\pp \, (cf.\ \cite[IV, Proposition 3.4]{Hu} and \cite[Appendix B]{CP1}).
Consider $\ga \colon A \to G/H$ and $\gb \colon G/H \to X$ as in
Definition \ref{def:gcat} and let $H \colon A \times I \to X$ be a
$G$-homotopy between $i \colon A \hookrightarrow X$ and $\gb 
\ga$. Since $G/H$ is a $G$-ANR, there exists a $G$-map $\gg \colon
W \to G/H$ where $W$ is a $G$-neighborhood of $A$ and $\gg|_A=\ga$.
Since $X$ is normal, we find a $G$-neighborhood $V$ of $A$ with
$\overline V\subset W$. We convert $X \times I$ to a $G$-space by
setting $g(x,t)=(g(x),t)$. Set
$$
P \,=\, X \times \{0\} \,\cup \, A \times I \,\cup\, \overline V \times \{1\}
\,\subset \, X \times I
$$
and define
$\gf \colon P \to X$ by
\[
\begin{array}{lll}
\gf (x,0) &=& x \,\text{ for } x \in X,\\
\gf (a,t) &=& H(a,t) \,\text{ for } (a,t) \in A \times I,\\
\gf(v,1)  &=& \gb ( \gg (v)) \,\text{ for } v \in \overline V.
\end{array}
\]
Then $\gf$ is well-defined and equivariant, and since $X \times
\{0\}$, $A \times I$ and $\overline V \times \{1\}$ are closed
subsets of $X\times I$, $\gf$ is continuous. Since $X$ is a $G$-ANR
and $P$ is a closed subspace of the normal $G$-space $X\times I$,
there exists a $G$-neighborhood $Q$ of $P$ in $X \times I$ and an
equivariant extension $\psi \colon Q \to X$ of $\gf$. For each $a\in
A$ the set $\{a \}\times I$ is compact, and so there is a
$G$-neighborhood $U_a$ of $a$ with $U_a \times I \subset Q$. Set $U =
\bigcup_a U_a \,\cap\, W$. Clearly, $U$ is a $G$-neighborhood of $A$
with $U\subset W$ and $U
\times I \subset Q$.

Now, $\psi|_{U \times I} \colon U \times I \to X$ yields a $G$-homotopy
between the inclusion $U \hookrightarrow X$ and $\gb \gg|_U$.
\proofend

\begin{proposition}  \label{p4:1}
If $A$ is a closed subspace of a binormal $G$-ANR $X$, then
\[
 \gcat A \,\ge\, \gccat A \,\ge\, \gccat_X A \,=\, \gcat_X A .
\]
\end{proposition}
\p
The first inequality follows from the simple fact that a closed subset
of a normal space is normal and from \text{Lemma \ref{l4:1}}, the
second inequality is clear since $A$ is closed, and the equality
follows from Lemmata \ref{l4:1} and \ref{l4:2}.
\proofend

\begin{remarks}  \label{r4}
{
\rm
{\bf 1.}
The closedness condition on $A$ in Lemma \ref{l4:1} is essential:
If $X$ is a circle, $x$ is a
point in $X$ and $A = X \setminus \{x\}$, then $2 = \ccat_X A > \cat_X
\! A =1$.

\s
{\bf 2.}
If $X$ fails to be an ANR, Proposition \ref{p4:1} might not hold:
If
\[
X = \bigcup_{n=1}^\infty \left\{ (x,y) \in \RR^2 \,\,\big|\,
\left(x- \textstyle{\frac{1}{n}} \right)^2 + y^2
= \textstyle{\frac{1}{n^2}} \right\}
\]
and $A = (0,0)$, then $\cat A =1$ while $\cat_X A = \infty$.
}
\end{remarks}

\section{Equivariant categories as index functions}

\ni
A $G$-map $\gf \colon X_1 \to X_2$ between $G$-spaces is a 
{\it $G$-homotopy equivalence} if there
exists a $G$-map $\psi \colon X_2 \to X_1$ such that 
$\psi \gf$ and $\gf \psi$ 
are $G$-homotopic to the identities.

\begin{lemma} \label{super}
Let $\gf \colon X \to X$ be a $G$-homotopy equivalence. If $U$ is an open
categorical subspace of $X$, then so is $\gf^{-1}(U)$.
\end{lemma}

\p Consider the commutative square
$$
\CD
\gf^{-1}(U) @>\gf>>U\\
@ViVV @VVjV\\
X@>\gf >> X
\endCD
$$
where $i$ and $j$ are the inclusions. Let $\psi \colon X \to X$ be a
$G$-homotopy inverse of $\gf$. If the composition
\begin{align*}
\!\!\!\!\!\!\!\!\!\!\!\!
\CD
U @>\ga>>  &  G/H @>\gb>> X
\endCD
\\
\intertext{is $G$-homotopic to $j \colon U \rightarrow X$, then the composition}
\!\!\!\!\!\!\!\!\!\!\!\!
\CD
\gf^{-1} (U) @>\ga \gf>>  &  G/H @>\psi \gb>> X
\endCD
\end{align*}
is $G$-homotopic to $i$, because 
$\psi \gb \ga \gf \simeq_G \psi j \gf = \psi \gf i \simeq_G i$.
\proofend

For any subset $A$ of a $G$-space $X$ we set $GA = \{ \,ga
\mid g \in G, \, a \in A \, \}$.

\begin{lemma}  \label{l6:1}
Let $X$ be a $G$-space and let $\gf \colon X \to X$ be a $G$-map.
Then, for every $N \in \NN$, each of the functions $2^X \times 2^X \to \NN\cup \{0\}$,
\[
\begin{array}{lll}
\nu_N^1 (A,Y) &=& \min \{ \gcat_X GA, N\},  \\       [.3 em]
\nu_N^2 (A,Y) &=& \min \{ \gcat_X (GA,GY), N\},  \\  [.3 em]
\nu_N^3 (A,Y) &=& \min \{ \gcat_X (GA \mmod GY), N\}  \\
\end{array}
\]
is an index function. Furthermore, the following holds.

\s
{\rm (i)} If $\gf$ is a $G$-homotopy equivalence, then $\nu^1_N$ is
$(\gf,Z)$-supervariant for every $Z \subset X$.

\s
{\rm (ii)}  If $\gf$ is $G$-homotopic to the identity, then $\nu^2_N$ is $(\gf,Z)$-supervariant for every $Z \subset X$.

\s
{\rm (iii)}
If there is a $G$-homotopy $\Phi_t$ between the
identity and $\gf$ with $\Phi_t(Z) \subset Z$ for all $t$,
then $\nu_N^3$ is $(\gf,Z)$-supervariant.
\end{lemma}

\p
The first claim is readily verified.

\s
(i) Let $\{A_1, \ldots, A_k\}$ be a covering of $G\gf(A)$ by open
categorical subspaces. By Lemma \ref{super},  $\{ \gf^{-1} (A_1), \ldots,
\gf^{-1} (A_k) \}$ is then a covering of $\gf^{-1}(G\gf(A))$ by open categorical subspaces. But $\gf^{-1}(G\gf(A))\supset GA$.

\s
(ii) Assume that $G\gf (A) \subset A_0
\cup A_1 \cup \dots \cup A_k$ where $A_1, \dots, A_k$ are open
$\cg$-categorical subspaces and $A_0$ is an open subspace which is 
$G$-deformable to $GZ$. Since $\gf$ is $G$-homotopic to the identity,
$\gf^{-1} (A_0)$ is $G$-deformable to $A_0$ and hence to $GZ$. 
We conclude now as in (i). 

\s
(iii)  Assume that $G \gf (A) \subset A_0
\cup A_1 \cup \dots \cup A_k$ where $A_1, \dots, A_k$ are open
$\cg$-categorical subspaces and $A_0$ is an open subspace which
contains $G \gf (A) \cap GZ$ and is $G$-deformable to $GZ
\mmod GZ$. Set $A_i' = \gf^{-1} (A_i)$. Then $GA \subset A_0' \cup A_1'
\cup \dots \cup A_k'$, and $A_1', \dots, A_k'$ are open
$\cg$-categorical subspaces. Let $\Phi_t$ be a $G$-homotopy between
the identity and $\gf$ with $\Phi_t (Z) \subset Z$. Then $GZ \subset
\gf^{-1}(GZ)$, and so $A_0'$ is an open subspace containing $GA \cap GZ$.
Moreover, the composition of $\Phi_t$ with a $G$-homotopy deforming
$A_0$ to $GZ \mmod GZ$ yields a $G$-homotopy deforming $A_0'$ to $GZ
\mmod GZ$.
\proofend

\section{LS theory for homotopy equivalences}

\ni
We let $G$ always act trivially on $\RR$. Thus, a $G$-function $f \colon 
X \ra \RR$ is a $G$-invariant function.
The {\it orbit type}\, of an orbit $Gx = G \{ x \}$ is its
$G$-homeomorphism type. Given a $G$-function $f \colon X \ra \RR$, two
orbits in $X$ are {\it equivalent}\, if $f$ has the same value on them
and if they have the same type and are $G$-deformable into each other
\cite[2.9]{CP1}. 
If $x$ is a fixed point of the $G$-map $\gf \colon X \ra X$, then the
whole orbit $Gx$ is fixed by $\gf$.

Recall that $\infty \ge \infty$ and $\infty \ge n$ for every $n \in
\ZZ$. We also agree that $\infty \ge \infty - n$ for every $n \in \ZZ$.

\begin{theorem}  \label{t5:1}
Consider a $G$-space $X$ and a class $\cg$ of homogeneous $G$-spaces.
Let $\gf \colon X \to X$ be a $G$-homotopy equivalence, and let $F$ be
the set of fixed points of $\gf$. Let  $- \infty < a < b<\infty$, 
and let $f \colon X \ra \RR$ be a $G$-function such that $(f, \gf)$ 
satisfies condition {\rm (D)} on $f^{-1}[a, b]$ and such that 
$\gcat_X f^a$ is finite. Set $F_d = F \cap f^{-1} (d)$.

\s
\begin{itemize}
\item[(a)] 
We always have
\[
\sum_{d \, \in \, ]a,b]} \gcat_X F_d \, \ge \, \gcat_X f^b - \gcat_X
f^a.
\]

\item[(b)]
If $X$ is a binormal $G$-ANR and $\cg$ contains all orbit types in $F
\cap f^{-1}]a,b]$, then $F \cap f^{-1}]a,b]$ contains infinitely many
$G$-orbits, or the number of equivalence classes of $G$-orbits in
$F\cap f^{-1}]a,b]$ is at least 
$\gcat_Xf^b - \gcat_Xf^a$.

\s
\item[(c)]
If $X$ is a binormal $G$-ANR and $f(F)\, \cap \, ]a,b]$ is discrete,
then
\[
\gcat F \cap f^{-1}]a,b] \, \ge \, \gcat_Xf^b - \gcat_Xf^a.
\]
\end{itemize}
\end{theorem}

\p
If $f(F) \, \cap \, ]a,b]$ is infinite, there is nothing to prove.
Suppose therefore that $f(F) \, \cap \, ]a,b] = \{ d_1,\dots,
d_m \}$ is finite.

\s
(a) If $\gcat_X f^b$ is finite, Lemma \ref{l6:1}\,(i) shows that
$$
\nu(A,Y):= \min \{ \gcat_X \! GA, \gcat_X f^b \}
$$
is a $(\gf,f^a)$-supervariant index function. Since $\nu(A,Y) =
\gcat_X \! GA$ for every
$A \subset f^b$, the claim then follows from Theorem \ref{t2}.
If $\gcat_X f^b = \infty$, we have to show that
\begin{equation}  \label{e:infty}
\sum_{i=1}^m \gcat_X F_{d_i} = \infty .
\end{equation}
Fix $N \ge \gcat_X f^a$. By Lemma \ref{l6:1}\,(i), $\nu_N (A,Y) = \min \{
\gcat_X \! GA, N \}$ is a $(\gf,f^a)$-supervariant index function. We have
$\nu_N (f^b, f^a) = N$
and $\gcat_X \! A \ge \nu_N (A)$ for every $A \subset X$. In view of
Theorem \ref{t2} we therefore conclude that
\begin{eqnarray*}
\sum_{i=1}^m \gcat_X F_{d_i} & \ge &
\sum_{i=1}^m \nu_N (F_{d_i}) \\
  &\ge& \nu_N (f^b,f^a) - \nu_N (f^a,f^a) \, = \,  N - \gcat_X f^a.
\end{eqnarray*}
Since $N \ge \gcat_X f^a$ was arbitrary, \eqref{e:infty} follows.

\s
(b) Assume that all the sets $F_{d_i}$ contain only finitely many
orbits. Let $E_1, \dots, E_{e_i}$ be the equivalence classes of orbits
in $F_{d_i}$. By assumption, all orbits in $E_j$ are $G$-deformable to
one of its elements, and $\cg$ contains this element. 
Lemma \ref{l4:2} thus implies that
$\gcat_X E_j =1$ for $j =1, \dots, e_i$. Therefore, $e_i =
\sum_{j=1}^{e_i}
\gcat_X E_j \ge \gcat_X F_{d_i}$.

\s
(c)
By Proposition \ref{p4:1},
\[
\gcat F \cap f^{-1}]a,b] \,=\, \sum_{i=1}^m \gcat F_{d_i} \,\ge\, \sum_{i=1}^m
\gcat_X F_{d_i} .
\]
\proofend

\begin{corollary}[{\rm cf.\ \cite[2.5]{CP1} and \cite[3.8 (1)]{CP2}}]
\label{c6:2}
Suppose that under the assumptions of Theorem {\rm\ref{t5:1}\,(b)} the set 
$f(F)\, \cap \, ]a,b]$ contains less than $\gcat f^b - \gcat f^a$ elements. 
Then one of the sets $F_d$ in $f^{-1} ]a,b]$ is not $G$-deformable to a $G$-orbit in $f^{-1} ]a,b]$.
\proofend
\end{corollary}

\begin{remarks}  \label{rs5:1}
{\rm
{\bf 1.} 
Proposition \ref{p4:1} implies that for a binormal
$G$-ANR, Theorem \ref{t5:1} holds with $\gcat$ replaced by $\gccat$.

\s
{\bf 2.}
We discuss the assumptions in Theorem \ref{t5:1}.

\s
(i)
The condition that $\gf$ is a $G$-homotopy equivalence 
cannot be omitted: 
If $X = S^1=\{x^2+y^2=1\}$, 
$\gf (X) = (0,-1)$ and $f(x,y) = y$ is the height function, then
$\sum_{d \,\in\, \RR} \cat_X F_d =1 < 2 = \cat X$.

\s
(ii) 
The condition that $b$ is finite cannot be omitted: 
Define $X$ to be the telescope (homotopy direct limit) of the sequence
$$
\CD
\cdots @>d>> S^1 @>d>>S^1 @>d>> S^1 @>d>> \cdots\,,\quad \deg d=2.
\endCD
$$
In greater detail, $X$ is the quotient space of the disjoint union
$$
Y=\coprod_{k=-\infty}^\infty S^1\times [k,k+1]
$$
under the following equivalence relation: $(z,k)\in S^1\times[k-1,k]$
is equivalent to $(z^2,k)\in S^1\times[k,k+1]$.

\s
We denote by $[z,k]\in X$ the image of $(z,k)\in S^1\times
[k-1,k]$ and by $[z,t]\in X$ the image of $(z,t)\in Y$ for $t\notin
\ZZ$. Consider the map
$$
\gf \colon X\to X,\quad  \gf[z,t]=[z, t-1]
$$
and the function $f \colon X \to \RR$, $f[z,t]=t$. Then $(f, \gf)$ clearly
satisfies condition (D), and it is easy to see that $\gf$ is a
homotopy equivalence. (This can be seen  directly or by observing that $X$ is
the Eilenberg--Mac Lane space $K(\ZZ[1/2],1)$ and that $\gf$ induces an
isomorphism of fundamental groups.) 
Furthermore, $F = \emptyset$.

Notice that $f^0$ is homotopy equivalent to $S^1$, and so $\cat_X f^0 =2$.
On the other hand, $\cat X > 2$; indeed, 
$\pi_1(X) = \ZZ[1/2]$ is not a free group, and
it is well known that the fundamental group of a
space whose category is at most $2$ is free \cite{F}. 
Therefore, $0= \sum_{d \,\in\, [0, \infty[} \cat_X F_d < \cat X - \cat_X
f^0$.

\s
(iii) 
Condition (D) cannot be omitted:
If $X = \,]0,1[$, $\gf (x) = \frac{x}{2}$ and $f(x) = x$, then $F =
\emptyset$ but $\cat X =1$.

\s
(iv) 
The discreteness condition on $f(F) \,\cap\, ]a,b]$ in (c)
cannot be omitted: If $X = S^1$ and 
$\gf \colon X \ra X$, $(x,y) \mapsto (\gf_1 (x,y), \gf_2 (x,y))$ is
such that $\gf (x,y) = (x,y)$ for $x \ge 0$ and $\gf_2(x,y) <y$ for
$x<0$, then $\gf$ is homotopic to the identity, and $f(x,y) =y$ is
a Lyapunov function for $\gf$. But $F = \{(x,y) \in X \mid x \ge 0\}$,
and so $\cat F =1 < 2 = \cat X$.
}
\end{remarks}

\section{LS theory for maps homotopic to the identity}

\ni
In this section we show that for maps homotopic to the identity, the
estimates in Theorem \ref{t5:1} can often be strengthened,
and that all estimates also hold for $b = \infty$.

\begin{theorem}  \label{t6:1}
Consider a $G$-space $X$ and a class $\cg$ of homogeneous $G$-spaces. 
Let $\gf \colon X \to X$ be a $G$-map which is
$G$-homotopic to the identity, and let $F$ be the set of fixed points of $\gf$. Let  $- \infty < a<b \le \infty$, and let $f \colon X \ra \RR$ be a $G$-function such that $(f, \gf)$ satisfies condition {\rm
(D)} on $f^{-1}[a, b]$. Set $F_d = F \cap f^{-1} (d)$.  

\m
{\rm (I)}
If $\gcat_X f^a$ is finite, then the statements 
{\rm (a), (b), (c)} in Theorem \ref{t5:1} also hold for $b = \infty$.

\s
{\rm (II)}
In case that $a \in f(F)$ and $a$ is an isolated point in $f(F) \cap [a,b[$, assume that a 
$G$-neighborhood of $f^a$ is $G$-deformable to $f^a$.
Then {\rm (a), (b), (c)} in Theorem \ref{t5:1} hold with 
$\gcat_X f^b - \gcat_X f^a$ replaced by $\gcat_X (f^b , f^a)$.

\s
\s
{\rm (III)}
If, in addition to the assumption in {\rm (II)},
there is a $G$-homotopy $\Phi_t$ between the identity and $\gf$ 
with $\Phi_t (f^a) \subset f^a$ for all $t$, then
{\rm  (a), (b), (c)} in Theorem \ref{t5:1} hold with 
$\gcat_X f^b - \gcat_X f^a$ replaced by $\gcat_X (f^b \mmod f^a)$.
\end{theorem}

\p
We may again assume that $f(F) \, \cap \, ]a,b] = \{ d_1,\dots, d_m \}$.
\begin{claim}  \label{c6:deformation} 
For any $c > d_m$ the space
$X$ is $G$-deformable to $f^c$. 
\end{claim}

\p  Choose a $G$-homotopy $\Phi \colon X \times [0,1] \ra X$ between the
identity and $\gf$.
For $x \in X$ and $t \ge 0$ set
\[
\widetilde \Phi (x,t) = \gf^{[t]} \left( \Phi (x, t - [t]) \right)
\]
where $[t] = \max \{ n \in \NN \cup \{0\}\mid n \le t \}$.
By Lemma \ref {l2:1}\,(i), for each $k \in \ZZ$ the number
\[
n_k = \min \{ n \in \NN \cup \{ 0 \} \mid \gf^n (f^k) \subset f^c \}
\]
is well-defined. Choose a non-decreasing continuous function $u \colon
\RR \ra [0, \infty[$ such that $u(k) = n_{k+1}$ for all $k \in \ZZ$.
Thus,
\begin{equation}\label{e:incl}
\gf^{[u(r)]}(f^r)\subset f^c
\end{equation}
for every $r \in \RR$. Define a $G$-function $h \colon X \ra
\RR_+$ by 
\[
h(x) = \max_{0 \le t \le 1} f( \Phi (x,t))
\]
and define $\Psi \colon X \times [0,1] \ra X$ by 
$
\Psi (x,t) = \widetilde\Phi (x, u(h(x)) \, t)$.
Then $\Psi$ is a $G$-homotopy between the identity and $\widetilde\Phi(x,
u(h(x)))$. We verify that $\widetilde\Phi (x, u(h(x))) \in f^c$ for each 
$x \in X$: Set $y = \Phi (x, u(h(x)) - [u(h(x))])$. Then $y \in
f^{h(x)}$. Therefore, by \eqref{e:incl},
$
\widetilde\Phi (x, u(h(x))) = \gf^{[u(h(x))]} (y) \in f^c.
$
\diam

\s
Claim \ref{c6:deformation}    
implies that 
$\gcat X = \gcat_X f^c$ and $\gcat_X (X,f^a) = \gcat_X (f^c,f^a)$.
Moreover, if the homotopy $\Psi$ above is constructed from a $G$-homotopy
$\Phi_t$ with $\Phi_t (f^a) \subset f^a$, then $\Psi (f^a,t) \subset
f^a$ for all $t$, and so $X$ is $G$-deformable to $f^c \mmod f^a$,
whence $\gcat_X (X \mmod f^a) = \gcat_X (f^c \mmod f^a)$.
We may therefore assume that $b$ is finite, and so, by the proof of
Theorem \ref{t5:1}, we are left with showing (a) in (II) and (III) for
$b$ finite.

\s
(II)(a) 
If $a \in f(F)$, a $G$-neighborhood of $f^a$ is $G$-deformable to
$f^a$ by assumption, and if $a \notin f(F)$, we conclude this by applying
\text{Lemma \ref{l2:1}\,(i)} to some $b \in \,]a, d_1[$ and to $U = \emptyset$.
So, $\gcat_X (f^a, f^a) = 0$. 
If $\gcat_X (f^b, f^a)$ is finite,  
we set 
\[
\nu (A,Y)   = \min \left\{ \gcat_X (GA, GY), \gcat_X (f^b, f^a) \right\}
\]
and conclude from Lemma \ref{l6:1}\,(ii) and Theorem \ref{t2} that 
$\sum_{i=1}^m  \gcat_X F_{d_i} \ge \gcat_X (f^b, f^a)$.
If $\gcat_X (f^b, f^a) = \infty$, we fix $N \in \NN$, set
$$
\nu_N (A,Y) = \min \left\{ \gcat_X (GA, GY), N \right\}
$$
and conclude from Lemma \ref{l6:1}\,(ii) and Theorem \ref{t2} that 
$\sum_{i=1}^m  \gcat_X F_{d_i} \ge N$.

\s
(III)(a) We use Lemma \ref{l6:1}\,(iii) to argue as for (II)(a).
\proofend

\begin{remarks}  \label{r6:1}
{\rm
{\bf 1.}
Theorem \ref{t6:1} shows that if $\gf$ is $G$-homotopic to the identity, 
then Corollary \ref{c6:2} can be refined and extended to the case $b = \infty$.

\s
{\bf 2.}
(a) and (c) of Theorem \ref{t5:1} and Theorem \ref{t6:1}\,(I) imply 
(a) and (c) of Theorem \ref{t0};
and since a connected and weakly locally contractible space is
path-connected, Theorem \ref{t0}\,(b) follows from Theorem \ref{t0}\,(a).

\s
{\bf 3.}
For a compact metric space $X$, Theorem \ref{t0} (a) with $\cat_X$
replaced by $\ccat_X$ and $\gf$ homotopic to the identity was stated by Ma\~n\'e \cite[Chapter II, Theorem 4.1]{M}.

\s
Proposition \ref{p4:1} implies that for a binormal $G$-ANR, 
Theorem \ref{t6:1}\,(I) holds with $\gcat$ replaced by $\gccat$.

\s
{\bf 4.}
We discuss the assumptions in Theorem \ref{t6:1}.

\s
(i)
Even if $b$ is finite, (II) in general does not hold for homotopy equivalences:
Let $S(k)$ be the unit circle in $\RR^2$ centred at $(0,2k)$
and set $X = \bigcup_{k \in \ZZ} S(k)$.
Define $\gf \colon X \ra X$ by $\gf (x,y) = (x, y-2)$ and set $f(x,y)
=y$. Then $F = \emptyset$ but $\cat_X (f^2, f^0) =1$.  

\s
(ii)
The assumptions in (II) and (III) cannot be omitted: Let 
\[
X = 
\left\{ (x,0) \mid -1 \le x \le 1 \right\}
\cup
\left\{ \left( \sin \textstyle{\frac{1}{y}}, y \right) \, \big| \, 
0 < y \le 1 \right\}  
\subset \, \RR^2 ,
\]
$\gf \left( \sin \textstyle{\frac{1}{y}}, y \right) = 
\left( \sin \textstyle{\frac{1 + 2 \pi y}{y}}, \frac{y}{1 + 2 \pi y} \right)$ 
and $\gf (x,0) = (x,0)$, and $f(x,y) = y$.
Then $F \cap f^{-1}]0,1] = \emptyset$ but $\cat_X (f^1, f^0) = \infty$.
Moreover, if
\[
X = 
\left\{ (x,y) \in S^1 \,\big|\, x \ge - \textstyle{\frac{\sqrt{3}}{2}},\, y \ge 0 \right\} ,
\]
$\gf (X) = (1,0)$ and $f(x,y) = y$, then $F \cap
f^{-1} ]\frac{1}{2},1] = \emptyset$ but $\cat_X \left(X \mmod
f^{\frac{1}{2}}\right) = 1$.

\s
{\bf 5.}
Replacing (i) in Definition \ref{def:gcat} by 

\s
\ni
\qquad \quad \;\,
$A_0$ is $G$-deformable to $Y$, \,and\, 
$A_0 \supset A \cap Y$

\s
\ni
we obtain a relative category $\gcat_X (A;Y)$ which is at least as large 
as $\gcat_X (A,Y)$. 
Since $\gf (f^a) \subset f^a$, the proof of Theorem \ref{t6:1}\,(II)
shows that if $b$ is  
finite, then the lower bound $\gcat_X (f^b,f^a)$ in (II) can be replaced 
by $\gcat_X (f^b;f^a)$.

\s
{\bf 6.}
(II) is not covered by (I) and (III): Let $X \subset \RR^3$ be the
space obtained by gluing the circle $X_1$ and the cylinder $X_2$
depicted below in the origin.
If $f$ is the height function  and $\gf$ restricts on $X_1$ to the
time-1-map of the negative gradient flow of $f$ and retracts $X_2$ onto
its bottom circle,  
then
$\cat_X (X,f^1) = 1 > 2-2 = \cat X - \cat_X f^1$,
and (III) does not apply -- indeed, $\sum_{d \,\in\, ]1,2]} \cat_X F_d =1 <2 = 
\cat_X (X \mmod f^1)$.

\begin{figure}[h] 
 \begin{center}
  \psfrag{x}{$x$}
  \psfrag{y}{$y$}
  \psfrag{z}{$z$}
  \psfrag{1}{$1$}
  \psfrag{X1}{$X_1$}
  \psfrag{X2}{$X_2$}
  \leavevmode\epsfbox{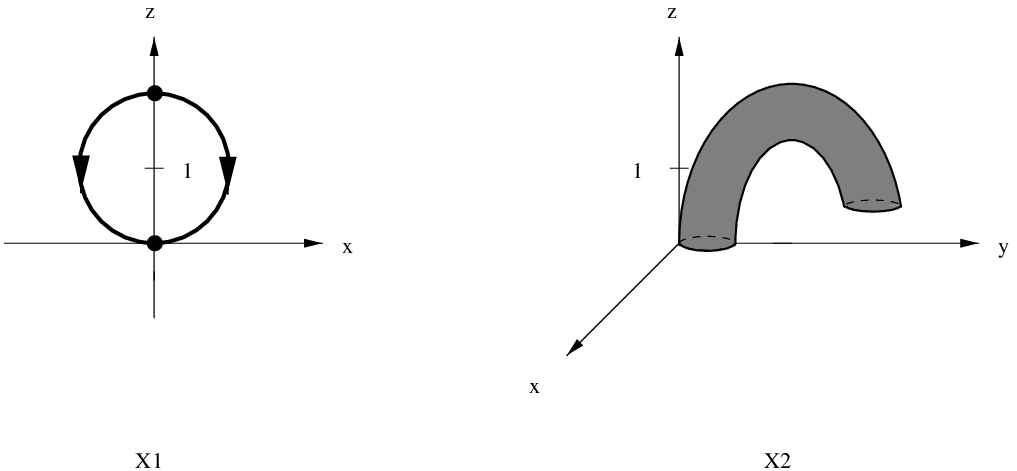}
 \end{center}
\end{figure}
%
%

(III) is not covered by (II): If $X = \{ (x,y) \in S^1 \mid y
\ge 0 \}$, 
$\gf \colon X \ra X$, $(x,y) \mapsto (\gf_1(x,y), \gf_2(x,y))$ 
is given by $\gf_2(x,y) = y^2$ and $x \, \gf_1(x,y) \ge 0$,  
and $f$ is the height function, then 
$\cat_X (X \mmod f^0) = 1 > 0 = \cat_X (X, f^0)$.

\s
{\bf 7.}
The main difficulty with (equivariant) categories is in their
computation. 
A lower bound for $\gcat_X (f^b \mmod f^a)$ is given by an equivariant
relative version of the cup-length, $\left( \cg, h^*
\right)\text{-cup-length}_X \left( f^b, f^a \right)$,
which depends on a choice of a $G$-equivariant cohomology theory $h^*$
\cite{Ba}.
Indeed, the sets $A_0, A_1, \dots, A_k$ used in the
definition of $\gcat_{\left( X,f^a \right)} f^b$ in \cite{CP2} are
required to be  
{\it relative open subsets}\, 
of $f^b$, whence
$$
\gcat_X \left( f^b \mmod f^a \right) \,\ge\, \gcat_{\left( X,f^a \right)} f^b,
$$
and the argument given in \cite[p.\ 58]{Ba} shows that
$$
\gcat_{\left( X, f^a \right)} f^b \,\ge\, \left( \cg, h^*
\right)\text{-cup-length}_X \left( f^b, f^a \right).
$$

A rougher but effectively computable lower bound for the category
$\gcat_X (f^b \mmod f^a)$
is the $\left( \cg, h^* \right)\text{-length}_X \left( f^b, f^a \right)$
\cite{CP2, Ba}.
}
\end{remarks}

\begin{definition}\label{def:semiflow}
\rm
A {\it $G$-semiflow} on a $G$-space $X$ is a family
$\Phi = \{\varphi_t\}$, $t \ge 0$, of $G$-maps $\varphi_t \colon X \to X$
such that $\gf_0 = id_X$ and $\varphi_s \circ \varphi_t =\varphi_{s+t}$
for all $s,t \ge 0$. Moreover, the map $X \times [0, \infty[ \, \to X$,
$(x,t) \mapsto \varphi_t(x)$ is required to be continuous.

A point $x \in X$ is called a {\it rest point} of $\Phi$ if
$\varphi_t(x)=x$ for all $t \ge 0$. 
The orbit of a rest point of $\Phi$ consists of rest points of $\Phi$
and is thus called a {\it rest orbit}\, of $\Phi$.

A $G$-semiflow $\Phi = \{\varphi_t\}$ is called {\it gradient-like}\, if there exists a $G$-function $f \colon X \to \RR$ such that
$f(\varphi_t(x)) < f(\varphi_s(x))$ whenever $t>s$ and $x$ is not a rest
point of $\Phi$. Any such function $f$ is called a {\it Lyapunov function} for $\Phi$.
\end{definition}

\begin{corollary}  \label{c6:flow}
Let $R$ be the set of rest points of the gradient-like $G$-semiflow
$\Phi = \{\gf_t\}$ on a $G$-space $X$. 
Let $f$ be a Lyapunov function for $\Phi$ which is
bounded below. 
Assume that there exists $\tau>0$ such that $(f, \gf_\tau)$ satisfies condition
{\rm (D)} on $X$.

\s
\begin{itemize}
\item[(a)]
We always have
$$
\sum_{d \, \in \, \RR} \gcat_X  R \cap f^{-1}(d) 
\, \ge \, \gcat X.
$$
\item[(b)]
If $X$ is a binormal $G$-ANR
and $\cg$ contains all orbit types in $R$, then $\Phi$ has at least
$\gcat X$ rest orbits.

\s
\item[(c)]
If $X$ is a binormal $G$-ANR and $f(R)$ is discrete, then
$$
\gcat R \,\ge\, \gcat X.
$$
\end{itemize}
\end{corollary}

\p
By the definition of a $G$-semiflow, $\gf_\tau$ is $G$-homotopic to the
identity. 
Denote the set of fixed points of $\gf_\tau$ by $F$.
\begin{claim}  \label{c:F=R}
 $F = R$.
\end{claim}
\p
The inclusion $R \subset F$ is obvious. Conversely, for $x \in F$ and $t 
\ge 0$, write $t = m \tau + r$ with $m \in \NN \cup \{ 0 \}$ and $r \in
[0, \tau[$. Then $\gf_{m \tau}(x) = \gf_\tau^m (x) = x$, and since
$\Phi$ is gradient-like, $\gf_r(x) = x$. Thus, 
$\gf_t(x) = \gf_r \left( \gf_{m \tau} (x) \right) = \gf_r (x) =x$.
\diam

\s
\ni
The corollary now follows from applying Theorem \ref{t6:1}\,(III) to
$(f, \gf_\tau)$. 
\proofend

We leave it to the reader to formulate a relative version of
the above corollary. 

\section{LS theory for homeomorphisms and flows}

\begin{definition}
\rm 
Fix a class $\BB$ of $G$-spaces.
Given a subset $A$ of a $G$-space $X$, we set
$
\BB_X A =  k
$
if $A$ can be covered by $k$ {\it open}\, subspaces $A_1,
\dots, A_k$ of $X$ such that every $A_i$ is $G$-homeomorphic to a
space from $\BB$, and if $k$ is minimal with this property. 
If there is no such number $k$, we set $\BB_X A = \infty$.
\end{definition}

\begin{theorem}  \label{t6:homeo}
Consider a $G$-space $X$. Let $\gf \colon X \to X$ be a $G$-homeomorphism,
and let $F$ be the set of fixed points of $\gf$. 
Let $-\infty < a < b < \infty$, 
and let $f\colon X \ra \RR$ be a $G$-function such that $(f,
\gf)$ satisfies condition {\rm (D)} on $f^{-1}[a, b]$ and such that
$\BB_X f^a$ is finite. 
Then
$$
\sum_{d \, \in \, ]a,b]} \BB_X \left( F \cap f^{-1} (d) \right) 
          \, \ge \, \BB_X f^b - \BB_X f^a .
$$
\end{theorem}

\p 
It is easy to see that for every $N \in \NN$ the function
\[
 2^X \times 2^X \to \NN \cup \{0\}, \quad \; (A,Z) \mapsto \min \{ \BB_X GA, N\}
\]
is an index function which is $(\gf,Z)$-supervariant for every
$Z\subset X$. The proof can now be completed as the one of Theorem
\ref{t5:1}\,(a).
\qed

\m
\begin{definition} 
\rm 
A semiflow $\{\gf_t\}$ on a space $X$ is called a
{\it flow}\, if each map $\gf_t$ is a homeomorphism.
\end{definition}

\begin{corollary}  \label{c8:1}
Let $R$ be the set of rest points of the gradient-like $G$-flow $\Phi =
\{\gf_t\}$ on a $G$-space $X$. Let $-\infty < a < b < \infty$, 
and let $f$ be a Lyapunov function for $\Phi$ such that there exists $\tau >0$
for which $(f, \gf_\tau)$ satisfies condition {\rm (D)} on $f^{-1}[a,b]$
and such that $\BB_X f^a$ is finite. 
Then
$$
\sum_{d \, \in \, ]a,b]} \BB_X \left( R \cap f^{-1} (d) \right) 
          \, \ge \, \BB_X f^b - \BB_X f^a .
$$
\end{corollary}

\p 
%
The corollary can be deduced from Theorem \ref{t6:homeo} in the same way
as Corollary \ref{c6:flow}\,(a) was deduced from Theorem
\ref{t6:1}\,(III). 
\proofend

\begin{remark}  \label{r8:a}
{\rm
Contrary to Corollary \ref{c6:flow}, the condition that $b$ is finite
cannot be omitted in Corollary \ref{c8:1}:
If $\BB$ is the class consisting of the annulus 
$\left\{ (x,y) \in \RR^2 \mid 1 < x^2 + y^2 < 4 \right\}$, and if
$X = \RR^2$ and $\Phi$ is the negative gradient flow of $f(x,y) = x^2 +
y^2$, then
$
\sum_{d \, \in \, ]1, \infty]} \BB_X \left( R \cap f^{-1}(d) \right) = 0 
< 2-1 = \BB_X X - \BB_X f^1$.
\diam
}
\end{remark}

We can play a similar game in the category of smooth manifolds.
We then consider a class $\BB$ of diffeomorphism types.
For example, let $B$ be the class consisting of $\RR^n$. 
Then $B_M A$ is the minimal number of open
balls necessary to cover $A \subset M$. 
Of course, $B_M A \ge \cat_M A$.
We set $B(M) = B_M M$, and
denote the minimal number of critical points of a smooth function on $M$ 
by $\Crit M$.
Following the proof of Theorem \ref{t6:homeo}, we obtain

\begin{corollary}  \label{c8:2}
Let $K$ be the set of critical points of a smooth
function on a closed smooth manifold $M$. Then
\[
\sum_{d \, \in ]a,b]} B_M \left( K \cap f^{-1}(d) \right) \, \ge \,
B_M f^b - B_M f^a .
\]
In particular, $\Crit M \ge B(M)$.
\proofend
\end{corollary}

We leave it to the reader to define a relative invariant $\BB_X ( A \mmod
Y)$ which can be used to refine Corollaries \ref{c8:1} and \ref{c8:2}.

\begin{remark}  \label{r8:1}
\rm
Singhof \cite{Sin} proved that $B(M) = \cat M$ for every closed smooth
$p$-connected manifold $M$ with 
$$
\cat M \, \ge \, \dfrac{n+p+4}{2(p+1)}
$$
provided $\cat M \ge 3$ and $\dim M \ge 4$. Moreover, it is easy to see
that if $B(M) = \cat M$, then $\cat(M \setminus \{x\})=\cat M-1$ for every
$x \in M$ \cite[p.\ 29]{Sin}. On the other hand, there is an example of a
closed manifold $Q$ with $\cat (Q\setminus \{x\})=\cat Q$, \cite{LSV}.
In particular, $B(Q) > \cat Q$. It is, however, still unknown
whether there are closed manifolds $M$ with $\Crit M > B(M)$.
\end{remark}

\section{The Palais--Smale condition (C) and condition (D)}

\ni
Consider a connected and complete Riemannian manifold $M$
without boundary modelled on a separable Hilbert space. 
For $m \in M$ we denote by
$\langle \,\, , \, \rangle_m$
the inner product on $T_m M$ and by $\| \; \|_m$
the norm induced by $\langle \,\, , \, \rangle_m$.
We say that a map is $C^{r,1}$, $r \ge 0$, if all its derivatives up to
order $r$ exist and are locally Lipschitz continuous. 
We assume that $M$ is of class $C^{1,1}$ and that the Riemannian metric
is $C^{0,1}$. Let $f \colon M \ra \RR$ be a $C^{1,1}$ 
function and let $Df(m)$ be its derivative at $m$. The $C^{0,1}$
vector field $\nabla f$, defined at $m$ by 
\[
D f (m) (v) = \langle \nabla f (m), v \rangle_m \quad \text{ for all } v \in
T_mM,
\]
is called the gradient vector field of $f$. 
A point $m \in M$ is called a critical point of $f$ if 
$\nabla f (m)$ vanishes.
The function $f$ is said to satisfy the Palais--Smale
condition (C) if the following holds:
\begin{quote}
If $A$ is a subset of $M$ on which $|f|$ is bounded but on which $\|
\nabla f \|$ is not bounded away from zero, then there is a critical
point of $f$ in the closure of $A$.
\end{quote}
\s
Choose a smooth, monotone function $g \colon \RR \ra \RR$
such that
\[
g(x) = 1 \,\text{ for } x \le 1, \quad g(x) \le x \,\text{ for } 1 \le x
\le
2, \quad g(x) = x  \,\text{ for } x \ge 2
\]
and set $h(m) = 1 / g( \| \nabla f(m) \|)$, $m \in M$.
Then $V = - h \, \nabla f$ is $C^{0,1}$
and bounded, and so, $M$ being complete, $V$
integrates to a flow $\{ \gf_t \}$, $t \in \RR$, on $M$. 

\begin{proposition}  \label{p:app}
Let $M$, $f$ and $\gf_\tau$, $\tau >0$, be as above.
If $f$ satisfies the Palais--Smale condition {\rm (C)},
then the pair $(f, \gf_\tau)$ satisfies condition {\rm (D)}.
\end{proposition}

\proof
Clearly, $f$ is a Lyapunov function for $\gf_\tau$.
Let $A \subset M$ and assume that $|f(a)| \le c < \infty$ for all $a \in
A$ and that
$\inf_{a \in A} \{f(a) - f(\gf_\tau (a))\} =0$.
We compute that for each $m \in M$,
\begin{eqnarray*}
\frac{d}{dt} f \left(\gf_t(m) \right) &=& Df \left( \gf_t(m) \right)
\left(
\frac{d}{dt} \gf_t(m) \right)    \\
&=& Df \left( \gf_t(m) \right) \left( V \left( \gf_t(m) \right)
\right) \\
&=& - h (\gf_t (m)) \left\langle \nabla f \left( \gf_t(m) \right),
\nabla f \left( \gf_t(m) \right) \right\rangle,
\end{eqnarray*}
and so
\begin{eqnarray}  \label{e:ps:1}
f (m) - f(\gf_\tau (m))  &=& - \int_0^\tau \frac{d}{dt} f \left( \gf_t(m)
\right) \, dt   \\
  &=& \int_0^\tau  h( \gf_t (m) ) \, \| \nabla f \left( \gf_t (m) \right)
\|^2
  \, dt . \notag
\end{eqnarray}
By assumption, there is a sequence $(a_n)_{n \ge 1} \subset A$ such
that
\begin{equation}  \label{e:ps:2}
f(a_n) - f(\gf_\tau (a_n))  <  \frac{\tau}{n} .
\end{equation}
Observe that
$ h(\gf_t(m)) \, \| \nabla f \left( \gf_t (m) \right) \|^2 <1$
only if
$\| \nabla f \left( \gf_t (m) \right) \| <1$.
Therefore,  \eqref{e:ps:1} and \eqref{e:ps:2} imply that there exists a
sequence
$(t_n)_{n \ge 1} \subset [0,\tau]$ such that
\begin{equation}  \label{e:ps:3}
\| \nabla f \left( \gf_{t_n} (a_n) \right) \|^2 < \frac{1}{n} .
\end{equation}
Set $b_n = \gf_{t_n} (a_n)$ and $B = \{b_n \}_{n \ge 1}$. Then, by
assumption and \eqref{e:ps:2},
\begin{equation}  \label{e:ps:4}
\begin{array}{lcl}
|f(b_n)| &\le& |f(a_n)| + |f(b_n) - f(a_n)| \\
         &\le& c+ |f(\gf_{t_n} (a_n)) - f(a_n)| \; < \; c+ \frac{\tau}{n}
\;
\le \;
          c+\tau .
\end{array}
\end{equation}
In view of \eqref{e:ps:3}, \eqref{e:ps:4} and condition (C),
we conclude that there exists $b^* \in \overline{B}$ with $\gf_\tau (b^*) =
b^*$. After passing to a subsequence, if necessary, we may assume that
$b^* = \lim_{n \ra \infty} b_n$.

For any $C^1$ path $\gs \colon [a,b] \ra M$ the length of $\gs$ is
defined by
$$
\int_a^b \left\| \frac{d}{dt} \gs (t) \right\| \, dt.
$$
For $m,m' \in M$ let
$d(m,m')$ be the infimum of the length of all $C^1$ paths joining $m$
and $m'$. The function $d$ thus defined is a metric on $M$ which is
consistent with the topology of $M$ \cite[\S 9]{P2}.
For each $n \ge 1$ the path $[0, t_n] \ra M$, $t \mapsto \gf_t
(a_n)$ is of class $C^1$. Therefore, by Schwartz's inequality, \eqref{e:ps:1} and
\eqref{e:ps:2},
\begin{eqnarray*}
d(b_n, a_n) \,=\, d(\gf_{t_n} (a_n), a_n) 
   &\le&  \int_0^{t_n} 
                      \left\| \frac{d}{dt} \gf_t (a_n) \right\| \, dt \\
   &\le&  \int_0^\tau \left\| \frac{d}{dt} \gf_t (a_n) \right\| \, dt \\
   &=& \int_0^\tau h (\gf_t(a_n)) \, \left\| \nabla f ( \gf_t (a_n) )
   \right\| \, dt \\ 
   &\le& \tau^\frac{1}{2} \left( \int_0^\tau  h (\gf_t(a_n))^2 \,
   \left\| \nabla f 
            \left( \gf_t (a_n) \right) \right\|^2 \, dt
\right)^\frac{1}{2}
\\
&\le& \tau^\frac{1}{2} \left( \int_0^\tau  h (\gf_t(a_n)) \,  \left\| \nabla f
            \left( \gf_t (a_n) \right) \right\|^2 \, dt
\right)^\frac{1}{2}
\\
   &\le& \frac{\tau}{\sqrt{n}} .
\end{eqnarray*}
We conclude that
\[
\lim_{n \ra \infty} d (b^* , a_n ) \le
\lim_{n \ra \infty} d (b^* , b_n ) + \lim_{n \ra \infty} d (b_n , a_n)
=0,
\]
i.e., $b^* \in \overline{A}$.
\proofend

A metric $G$-space $X$ is called a {\it metric $G$-ANR}\, if
for each closed $G$-subspace $Z$ of a {\it metric}\, $G$-space $Y$
every $G$-map $\rho \colon Z
\to X$ admits an equivariant extension to a $G$-neighborhood of $Z$.

\begin{corollary}  \label{c:sch}
Let $M$ be a Riemannian Hilbert manifold as above.
Suppose that $M$ is a $G$-space
and that the Riemannian metric on $M$ is $G$-invariant. 
Let $K$ be the set of critical points of a $C^{1,1}$ $G$-function
$f \colon M \ra \RR$ which satisfies condition {\rm (C)} on $M$ and is
bounded below.
\s
\begin{itemize}
\item[(a)]
We always have
$$
\sum_{d \, \in \, \RR} \gcat_M  K \cap f^{-1}(d) 
\, \ge \, \gcat M.
$$
\item[(b)]
If $M$ is a metric $G$-ANR and $\cg$ contains all orbit types in $K$, 
then $f$ has at least $\gcat M$ critical orbits.  

\s
\item[(c)]
If $M$ is a metric $G$-ANR and $f(K)$ is discrete, then
$$
\gcat K \,\ge\, \gcat M.
$$
\end{itemize}
\end{corollary}

\p
Going once more through Section 4, we see that Proposition \ref{p4:1}
holds for metric $G$-ANR's $X$. Therefore, Corollary \ref{c6:flow} holds for 
metric $G$-ANR's.
Since the Riemannian metric on $M$ is $G$-invariant, 
$\nabla f$ and $V$ are $G$-equivariant, and so is $\gf_t$, $t \in \RR$.
Moreover, the critical points of $f$ are the rest points of $\gf_1$. 
The corollary now follows in view of Proposition \ref{p:app}.
\proofend

Localizing Proposition \ref{p:app} and applying the relative version of 
Corollary \ref{c6:flow}, we obtain a relative version of the above corollary.

\begin{remarks}  \label{rs9:1}
\rm
{\bf 1.}
Let $M$ be as in the basic assumptions of Corollary \ref{c:sch}. 
Since $M$ is metrizable, it is paracompact.
Appendix B of \cite{CP2} thus provides sufficient conditions for $M$
being a metric $G$-ANR.
They are, e.g., fulfilled if $G$ acts trivially, or if $M$ and the action 
$G \times M \ra M$ are $C^{2,1}$.

\s
{\bf 2.}
By the above remark, Corollary \ref{c:sch}\,(b) recovers  
the basic Lusternik--Schnirelmann theorem
for Hilbert manifolds first obtained in \cite{Sch}.

\s
{\bf 3.}
Versions of (a) and (b) in Corollary \ref{c:sch} have been proved for
$C^{1,1}$ functions on complete $C^{1,1}$ Finsler manifolds 
(\cite[Theorems 7.1 and 7.2]{P4} and \cite[$\S$ \!3]{CP2}) and for continuous
functions on weakly locally contractible complete metric spaces 
(see \cite{CD}).
We notice that in these situations also the analogs of (c) hold.
\end{remarks}

\enddocument